\documentclass[12pt]{amsart}

\usepackage{latexsym, amssymb, amsmath}
\usepackage{amsmath}
 \usepackage{eucal}

    \newtheorem{rema}{Remark}[section]
    \newtheorem{propo}[rema]{Proposition}

   \newtheorem{theo}[rema]{Theorem}
 \newtheorem{conj}[rema]{Conjecture}
   \newtheorem{defi}[rema]{Definition}
    \newtheorem{lemma}[rema]{Lemma}
    \newtheorem{corol}[rema]{Corollary}

     \newcommand{\nno}{\nonumber}
	\newcommand{\lbar}{\big\vert}
        
	\newcommand{\p}{\partial}

 \newcommand{\pf}{{\it Proof:}\hspace{2ex}}
 
 \newcommand{\epfv}{\hspace{1em}$\Box$\vspace{1em}}

\newcommand{\bC}{{\mathbb C}}

\newcommand{\bN}{{\mathbb N}}

\newcommand{\BQ}{\begin{eqnarray}}
\newcommand{\EQ}{\end{eqnarray}}
\newcommand{\BQn}{\begin{eqnarray*}}
\newcommand{\EQn}{\end{eqnarray*}}

	\newcommand{\Pz}{\frac {\partial}{\partial z}}

\begin{document}
    \bibliographystyle{alpha}

\title[Exponential Formulas for Jacobians and Jacobian Matrices]
{Exponential Formulas for the Jacobians and Jacobian Matrices of 
         Analytic Maps}
      \author{Wenhua Zhao}      

\begin{abstract}
Let $F=(F_1, F_2, \cdots F_n)$ be an $n$-tuple of formal power series in $n$ variables
of the form $F(z)=z+ O(|z|^2)$. It is known that, 
there exists a 
unique formal differential operator
$A(z)=\sum_{i=1}^n a_i(z)\frac {\p}{\p z_i}$ such that
$F(z)=exp (A)z$ as formal series. In this article, 
we show the Jacobian ${\mathcal J}(F)$ and the Jacobian matrix $J(F)$ of $F$
can also be given by some exponential formulas. Namely,
${\mathcal J}(F)=\exp (A+\triangledown A)\cdot 1$, where 
$\triangledown A(z)=
\sum_{i=1}^n \frac {\p a_i}{\p z_i}(z)$, and 
$J(F)=\exp(A+R_{Ja})\cdot I_{n\times n}$, where $I_{n\times n}$ is the identity matrix
and $R_{Ja}$ is the multiplication operator by $Ja$ for the right.
As an immediate 
consequence, we get an elementary proof for the known 
result that ${\mathcal J}(F)\equiv 1$ 
if and only if $\triangledown A=0$.  Some consequences and
applications of the exponential 
formulas as well as their relations with the well known Jacobian Conjecture 
are also discussed.
\end{abstract}

\thanks{{\it 2000 Mathematics Subject Classification}. 
32H02, 32A05, 14R15}

\maketitle


\renewcommand{\theequation}{\thesection.\arabic{equation}}
\renewcommand{\therema}{\thesection.\arabic{rema}}
\setcounter{equation}{0}
\setcounter{rema}{0}

\renewcommand{\theequation}{\thesection.\arabic{equation}}
\renewcommand{\therema}{\thesection.\arabic{rema}}
\setcounter{equation}{0}
\setcounter{rema}{0}
\setcounter{section}{0}

\section{Introduction}\label{Sec-1}

This research work mainly motivated by the well known Jacobian Conjecture 
and inspired by an exponential formula in Conformal Field theory. First let us
recall 

{\bf Jacobian Conjecture:}\quad
 {\it Let $k$ be a field of characteristic $0$ and 
$F: k^n\to k^n$ be a polynomial map.
If Jacobian $j(F)=Det \left (\frac {\p F_i}{\p z_j} \right ) = 1$, 
then $F$ is an automorphisms 
whose inverse is also a polynomial map.}

This conjecture was first proposed by O. H. Keller in $1939$.
For the history of this conjecture, 
(See \cite{BCW}, \cite{W1} and \cite{M} and references there).
Since then it has been
attracting enormous efforts from mathematicians. 
But unfortunately, this conjecture 
remains widely open at the present time. 
Nevertheless, many important results have been
obtained in last six decades from the efforts of 
mathematicians trying to solve Jacobian
Conjecture. Some of these results are 
not only crucial to the Jacobian Conjecture, 
they also play very important roles in other 
mathematical research areas. 

One of the effective approaches to the Jacobian Conjecture is 
to develop nice formulas for the formal inverse $G$ of 
the polynomial map $F$ and to see 
if it is also a polynomial map. Several important formulas 
have been found and well studied, 
among which the most well known 
are Abhyankar's inversion formula 
(See \cite{Ab}) and the tree expansion
formula for the formal inverse $G$ of $F$. 
(See \cite{BCW} and \cite{W2}). 

Interestingly, an exponential formula 
for the formal power series or 
holomorphic functions in one
variable, which plays a crucial role 
in two dimensional Conformal Field Theory, 
seems closely related with 
the Jacobian Conjecture. To be more precise, 
let $F(x)=x+O(x^2)$ be a formal power series in one variable $x$. 
Then there exists a unique formal differential operator 
$A(x)=a(x)\frac {\p}{\p x}$ with $o(a)\geq 2$ 
such that $F(x)= e^A x$. (Note that the exponential formula we 
quote here is a little different from the one used in \cite{TUY}).
The main reason that the exponential formula above is so important in 
two dimensional Conformal Field Theory is that it 
gives the Virasoro algebra structure,  
which is the most fundamental algebraic structure to the whole theory. 
For more detail, see \cite{TUY}, \cite{H} and \cite{Z1}.

One of the advantages of the exponential formula $F(x)= e^A x$
for the formal power series $F(x)$
is that $e^A$ is an automorphism of
the algebra $\bC [[x]]$ of 
formal power series in one variable. 
This is because that 
$A$ itself is 
a derivation  of
the algebra $\bC [[x]]$ and it is well known in 
Lie algebra theory that the exponential of 
any derivation of an algebra is an automorphism of
the algebra. 
As an immediate consequence of this observation, 
the formal inverse $G$ of $F$ is given 
by the exponential formula
$G(x)=e^{-A} x$. 
Regarding the Jacobian Conjecture,
it is certainly very interesting to see that the formal inverse 
$G$ of $F$ is given in such a simple way. Actually, for the
formal power series in several variables, 
we also have similar exponential formulas 
(See \cite{P} and also Proposition \ref{6P2.1}). Namely,
let $F=(F_1, F_2, \cdots F_n)$ be an $n$-tuple of formal 
power series in $n$ variables
of the form $F(z)=z+ O(|z|^2)$. 
Let $G=(G_1, G_2, \cdots G_n)$ be the formal inverse of $F$, i.e.
$F(G)=G(F)=z$, where $z=(z_1, z_2, \cdots, z_n)$. Then
there exists a unique formal differential operator
$A=\sum_{i=1}^n a_i(z)\frac {\p}{\p z_i}$ with $o(a_i(z))\geq 2$ such that
$F_i(z)=\exp (A)z_i$ $(i=1, 2, \cdots, n)$. 
By the similar reason, $e^A$ is an automorphism of
the algebra $\bC [[z]]$ of 
formal power series in $z$ and
$G_i(z)=e^{-A} z_i$ for $i=1, 2, \cdots, n$.

Since the formal power series $F$ as well as 
its formal inverse $G$ are completely determined by
a unique formal differential operator $A$, 
naturally one may ask: how  
does the formal differential operator $A$ 
determine the Jacobian ${\mathcal J}(F)$ 
and Jacobian matrix $J(F)$ of $F$? or in other words, are there any formulas via which
the differential operator $A$ also completely determines ${\mathcal J}(F)$ and
$J(F)$?
In this article, we show that the answer to the question above is ``yes". 
In Section \ref{Sec-2}, we give two exponential formulas for 
the Jacobian ${\mathcal J}(F)$ 
and Jacobian matrix $J(F)$ of $F$, respectively. To be more precise, 
in Theorem \ref{MT}, we show that 
${\mathcal J}(F)=exp (A+\triangledown A)\cdot 1$, where 
$\triangledown A(z)=
\sum_{i=1}^n \frac {\p a_i}{\p z_i}(z)$ is 
the {\it divergence} of the operator $A$. In Theorem \ref{JM-MT},
We show that 
$J(F)=\exp(A+R_{Ja})\cdot I_{n\times n}$, where $I_{n\times n}$ is the identity matrix
and $R_{Ja}$ is the multiplication operator by $Ja$ for the right.
As an immediate 
consequence, we also give an elementary proof for the known 
result that ${\mathcal J}(F)\equiv 1$ 
if and only if $\triangledown A=0$. (See Corollary \ref{C2.10}). 
Various interesting properties of 
the differential operator $A$ and the formal 
deformation $F_t(z)=e^{tA}z$ are also derived in this section.

In Section \ref{Sec-3}, we first give some explanations 
about the exponential formulas derived in 
Theorem \ref{MT} and Theorem \ref{JM-MT} by relating theom with 
some well known formula in linear algebra. Then we study 
the consequences of these exponential formulas 
to the Jacobian Conjecture, especially, we give a new proof to 
a theorem of Bass, Connell and Wright, in \cite{BCW}.  (See Theorem \ref{Thm-BCW}).

In Section \ref{Sec-4}, we discuss some open 
problems related with these exponential formulas
and the Jacobian Conjecture.

 Most of the results of this article comes form
the third topic of the author's Ph.D Thesis \cite{Z1} in
the University of Chicago, except Theorem \ref{JM-MT}, Theorem \ref{Thm-BCW}
and the "Explanation" part of  
Section \ref{Sec-3}
are added later. Theorem \ref{MT} is also given 
in a more general form than the one in \cite{Z1}.
The author is very grateful to his Ph.D advisor, 
Professor Spencer Bloch for encouragement, discussions and  
pointing out an error in the early 
version of this work. The author is very thankful to
Professor Yi-Zhi Huang for many
personal communications and the suggestion to the author the last open problem in
Section \ref{Sec-4}. The author thanks Professor Xiaojun Huang for many
discussions on some 
analytic aspects of this work. Great appreciation also goes 
to the Department of Mathematics, the University of Chicago 
for financial supports during the author's graduate study.

\renewcommand{\theequation}{\thesection.\arabic{equation}}
\renewcommand{\therema}{\thesection.\arabic{rema}}
\setcounter{equation}{0}
\setcounter{rema}{0}

\section{Exponential Formulas}\label{Sec-2}

{\bf Notation:}

$(1)$ Let $z_1, z_2, \cdots z_n$ be $n$ commutative variables and 
$z=(z_1, \cdots z_n )$. Let 
$\bC [[z]]=\bC [z_1, z_2, \cdots z_n]$ be the 
algebra of polynomials in $n$ variables, 
$\bC [[z]]$ be the algebra of formal power series.
For any $k\geq 0$, set $\bC_k[[z]]=
\bC_k[[z_1, z_2, \cdots z_n]]$ be the subalgebra 
consisting of 
the elements of $\bC [[z]]$
whose lowest degree is greater or equal to $k$.
 
$(2)$ For any $F=(F_1, F_2, \cdots , F_n) \in 
\bC [[z]]^n$, set 
\BQ
JF(z)=\left (\frac {\partial F_i}{\partial z_j} \right )_{1\leq i, j \leq n}
\EQ
\BQ
{\mathcal J}F(z)=
Det \left (\frac {\partial F_i}{\partial z_j} \right )_{1\leq i, j \leq n}
\EQ

We call $JF$ the {\it Jacobian matrix}  and ${\mathcal J}F(z)$ the 
{\it Jacobian} of $F$.

Let ${\mathcal F}_1$ be the set of the elements 
$F=(F_1, F_2,\cdots, F_n)\in \bC [[z]]^n$ 
such that
$F_i(z)=z_i+\text {high degree terms}$, 
for $i=1, 2, \cdots n$. Note that for any
analytic map $F: U \to \bC^n$ with Jacobian ${\mathcal J}(F)(0)\neq 0$ for the
some open neighborhood $U$ of $0\in \bC^n$, composing with some line
isomorphism if necessary, the formal series of $F$ will be in ${\mathcal F}_1$. 
Another observation is that, any $F\in {\mathcal F}_1$ gives an automorphism
of the algebra $\bC [[z]]$, which sends $z_i$ to $F_i$. The inverse of this
automorphism is the automorphism induced by the formal inverse of $F$.

One remark is that all the proofs and results in this paper work equally well for
any field of characteristic $0$, not necessarily algebraic closed.
But for convenience, we will always take $\bC$ to be the ground field. 

The following proposition is known. For example, see \cite{P}. Here we give an
elementary proof.

\begin{propo}\label{6P2.1}
For any $F=(F_1, F_2, \cdots , F_n) \in {\mathcal F}_1$, there exists a unique 
$a=(a_1, a_2, \cdots , a_n) \in 
{\bC_2}[[z]]^n$ such that 
\BQ\label{EF1}
F_i(z)=exp({a(z)\frac d{dz}})z_i=exp(A)z_i
\EQ
where 
\BQ
A(z)&=&a(z)\frac d{dz}=\sum_{i=1}^n a_i(z)\frac {\p}{\p z_i}\\
exp (A)&=&exp ({a(z)\frac d{dz}})=\sum_{k=0}^\infty \frac {(a(z)\frac d{dz})^k}{k!}
\EQ
\end{propo}

\pf  
This can be checked directly by solving the formal equation (\ref{EF1}) 
incursively as following.

For $i=1, 2, \cdots , n$, we write 
\BQ
F_i(z)&=&z_i+b_i^{(2)}(z)+b_i^{(3)}(z)+\cdots +b_i^{(k)}(z)+\cdots \\
a_i(z)&=&a_i^{(2)}(z)+a_i^{(3)}(z)+\cdots +a_i^{(k)}(z)+\cdots
\EQ
where $a_i^{(k)}(z)$ and $b_i^{(k)}(z)$, for any $k\in \bN$,  are homogeneous 
polynomials of degree $k$. We also write
$F^{(k)}=(F_1^{(k)}, F_2^{(k)} \cdots ,F_n^{(k)})$,     
$a^{(k)}=(a_1^{(k)}, a_2^{(k)} \cdots ,a_n^{(k)})$ and
$A^{(k)}=a^{(k)}\frac {\p}{\p z}=\sum_{i=1}^n a_i^{(k)}\frac {\p}{\p z_i}$.
Notice that the operator $A^{(k)}$ increase degree by $k-1$.

From the equations (\ref{EF1}), we get 
\BQ\label{PE1}
&{}& z_i+\sum_{k=1}^\infty \frac {(a(z)\frac d{dz})^k}{k!}z_i \\
&{}& \quad \quad 
=
z_i+b_i^{(2)}(z)+b_i^{(3)}(z)+\cdots +b_i^{(k)}(z)+\cdots \nno
\EQ
Comparing the homogeneous parts of both sides of (\ref{PE1}), we get 
\BQ
a_i^{(2)}&=& b_i^{(2)} \nno \\
a_i^{(3)}&=& b_i^{(3)}- 
\sum_{k=1}^n a_k^{(2)} \frac {\p a_i^{(2)}}{\p z_k} \nno\\
&\cdots &  \nno \\
a_i^{(m)}&=& b_i^{(m)}- \sum_{1\leq r<m}
\sum_{\begin{subarray}{1}
k_1+k_2+\cdots k_r=m+r-1 \\ \label{Incursive}
 k_1, k_2, \cdots k_r\geq 2 \end{subarray} }
\frac {A^{(k_1)}A^{(k_2)}\cdots A^{(k_r)}}{k_1! k_2! \cdots k_r!}z_i
\EQ
Hence $a(z)$ is completely determined by the equations above.
\epfv

One easy corollary of the calculation above is the following

\begin{corol}\label{L-odd}
$F$ is odd if and only if $a(z)$ is odd.
\end{corol}

This can also be proved by the similar arguments for 
Proposition \ref{P-even}.

\begin{defi}
We call the formal differential operator $A$ in Proposition \ref{6P2.1}
the {\it associated differential operator}
of $F$. We also define 
\BQ
(\triangledown A)=(\triangledown a)(z)=\sum_{i=1}^n \frac {\p a_i}{\p z_i}(z)
\EQ
and call it the {\it divergence} of the differential operator $A$.
\end{defi}

One of the advantages of the formula (\ref{EF1}) is that 
the operator $exp(A)$ or $exp ({a(z)\Pz})$ is an automorphism of the $\bC$-algebra
$\bC [[z]]$ which maps $z_i$ to $F_i$. This follows from the well known fact 
that the exponential of any derivative of any algebra, when it is well defined,
is an automorphism of that algebra. It is because this remarkable 
property that the formula (\ref{EF1}) in the case of one variable 
plays a very important role in conformal field
theory. See \cite{H} and \cite{TUY}. (The formula used in \cite{TUY} is a little
different from (\ref{EF1})).
The following are some immediate
consequences of the property above.

\begin{lemma}\label{L2.5}
$a)$ Let $F^{-1}=(F_1^{-1}, F_2^{-1}, \cdots ,  F_n^{-1})$ be 
the formal inverse of $F$, i.e. the composition 
$F\circ F^{-1}=F^{-1}\circ F$ 
is identity map of $\bC [[z]]$.
Then 
\BQ
F^{-1}(z)=exp ({-A(z)})z=exp ({-a(z)\Pz}) z
\EQ

$b)$ For any element $g(z)\in \bC [[z]]$, we have
\BQ
g(F(z))=exp ({a(z)\Pz}) g(z)
\EQ
In particular, for any $k\geq 0$, we have
\BQ
 F^{[k]}(z)=exp ({kA(z)})z=exp ({ka(z)\Pz}) z
\EQ
where
\BQ
F^{[k]}(z)=\underbrace {F\circ F\circ \cdots F}_{k \, \,  \text{copies} }
\EQ
is the $k^{th}$-power of the automorphism of $\bC [[z]]$
defined by $F$ which sends
$z_i$ to $F_i$.
\end{lemma}

Another advantage of the formula (\ref{EF1}) is that it allows us 
to deform the formal power series $F$ in a very natural way. Introduce 
another variable $t$ which commutes with $z_i$ 
and define 
\BQn
F_t(z)=F(z; t)=(F_1(z; t), F_2(z; t), \cdots ,F_n(z; t))
\EQn
by setting
\BQ\label{Ft}
F_i(z; t)= exp ({tA(z)}) z_i=exp ({ta(z)\Pz})z_i
\EQ
Note that $F_i(z; t)\in \bC [t][[z]]$, i.e. it is a formal power series in 
$\{ z_i\}$ with coefficients in $\bC [t]$. In particular,
for any $t_0\in \bC$, $F(z; t_0)\in {\mathcal F}_1$ and 
when $t=k\in \bN$, 
$F(z; k)$ is just the $k^{th}$-power $F^{[k]}$ 
of the isomorphism $F$.
This deformation will play the key role in our
later arguments.

\begin{lemma}\label{L2.6}
For any $g(z; t)\in \bC [t][[z]]$, 
\BQ\label{g(z; t)}
\frac {\p}{\p t}g(z; t)=Ag(z; t)
\EQ
if and only if $g(z; t)=u(F(z; t))=exp (tA)u(z)$ for some $u\in \bC [[z]]$.
\end{lemma}

\pf  First let $g(z; t)=exp (tA)u(z)$, then
\BQn
\frac {\p}{\p t}g(z; t)&=&\frac {\p}{\p t}exp (tA)u(z) \\
&=& Aexp (tA)u(z)\\
&=& Ag(z; t)
\EQn

Conversely, suppose that $g(z; t)$ satisfies (\ref{g(z; t)}).  then, 
 by chain rule, we have 
\BQn
\frac {\p}{\p t}exp(-tA)g(z; t)&=&
-Aexp(-tA)g(z; t)+exp(-tA)\frac {\p}{\p t}g(z; t) \\
&=&-Aexp(-tA)g(z; t)+Aexp(-tA)g(z; t) \\
&=& 0
\EQn
So $exp(-tA)g(z; t)$ does not depend on $t$ and is in $\bC [[z]]$.
Set 
$u(z)=exp(-tA)g(z; t)$, we have $g(z; t)=exp(tA)u(z)$.
\epfv

The following property is a little bit strange.

\begin{propo}
\BQ \label{PE9}
J(F)(z; t) 
\begin{pmatrix} 
 {a_1(z)}\\
{a_2(z)}\\
\vdots \\
 {a_n(z)}
\end{pmatrix} =
\begin{pmatrix} 
{a_1(F(z; t))}\\
{a_2(F(z; t))}\\
\vdots \\
 {a_n(F(z; t))}
\end{pmatrix}
\EQ
or in short notations
\BQ
AF(z; t)=J(F)(z; t)a(z)=a(F(z; t))
\EQ
\end{propo}

\pf  This follows from the following straightforward calculations.
Consider
\BQ
\frac {\p}{\p t}F_i(z; t)&=& \frac {\p}{\p t}exp (ta(z)\Pz)z_i \nno\\
&=&\sum_{k=1}^n a_k(z)\frac {\p}{\p z_k}exp (ta(z)\Pz)z_i\nno \\ 
&=&\sum_{k=1}^n a_k(z)\frac {\p}{\p z_k}F_i(z; t) \nno \\ \label{PE7}
&=&\sum_{k=1}^n\frac {\p F_i(z; t)}{\p z_k}a_k(z)
\EQ
On the other hand, note that the operators $a(z)\Pz$ and $exp(a(z)\Pz)$ commute 
with each other, so we also have
\BQ
\frac {\p}{\p t}F_i(z; t)&=&exp (ta(z)\Pz)(\sum_{k=1}^n a_k(z)\frac{\p}{\p z_k})z_i\nno \\
&=&exp (ta(z)\Pz)a_i(z) \nno \\
&=&a_i(exp (ta(z)\Pz)z) \nno \\ \label{PE8}
&=& a_i(F(z; t))
\EQ

Comparing (\ref{PE7}) and (\ref{PE8}), we get (\ref{PE9}).
\epfv

Unfortunately, the equation (\ref{PE9}) does not completely determine the operator
$A(z)=a(z)\Pz$. Instead we have the following explicit formulas for
$a(z)$ and the inverse $G=(G_1, G_2, \cdots , G_n)$ of $F$.

\begin{propo}\label{Explicit}
$a)$ 
\BQ
a(z)=-\sum_{k=1}^\infty \frac 1k (1-e^A)^kz=
-\sum_{k=1}^\infty \frac 1k \left ( 
\sum_{j=0}^k (-1)^j \binom kj F^{[j]}(z) \right )
\EQ
$b)$
\BQ\label{GE}
G(z)=z+\sum_{k=1}^\infty (1-e^A)^kz=z+\sum_{k=1}^\infty 
\left ( 
\sum_{j=0}^k (-1)^j \binom kj F^{[j]}(z) \right )
\EQ 
\end{propo}

Notice that the operator $1-e^A$ strictly increases the degree, so the infinite sums
that appear in the lemma above all make sense.

\pf  $a)$ follows from the following formal identity
\BQ
A=log \, e^A=log (1-(1-e^A))=-\sum_{k=1}^\infty \frac 1k (1-e^A)^k
\EQ

$b)$ Since the formal inverse of $F$ exists and is unique, it is enough to check
that the formal series $G$ given by  (\ref{GE}) is the inverse of $F$.

Consider 
\BQn
(G\circ F)(z)&=&e^AG(z) \\
&=&e^Az+ \sum_{k=1}^\infty (1-e^A)^ke^Az \\
&=&e^Az +\sum_{k=1}^\infty (1-e^A)^{k}z -\sum_{k=1}^\infty (1-e^A)^{k+1}z\\
&=&e^Az +(1-e^A)z\\
&=&z
\EQn
\epfv

Now we begin to prove our exponential formula for the Jacobian ${\mathcal J}(F_t)$.

\begin{theo}\label{MT}
$(a)$ In the notations above, we have
\BQ\label{MF}
{\mathcal J}(F_t)(z)=exp (t{a(z)\frac d{dz}+ t\triangledown a(z)}) \cdot 1 
\EQ
where $F_t(z)=F(z; t)=(F_1(z; t), F_2(z; t), \cdots , F_n(z; t))$ as before.

$(b)$ For any $u\in \bC [[z]]$, we have
\BQ \label{G-MF}
exp (tA + t\triangledown a(z))u
=u(F(t, z)){\mathcal J}F(t, z) 
\EQ
\end{theo}

It is easy to see that $(a)$ is an immediate consequence of $(b)$, but here 
we need prove $(a)$ first.

\pf $(a)$ 
To keep notations simple, here we only give the proof for
the case of two variables. For the general
cases, the ideas are completely same. 

Let  $K(t)=
exp ({ta(z)\frac d{dz}+ t\triangledown a(z)})\cdot 1 $
and $H(t)={\mathcal J}(F_t)$, i.e. the Jacobian of $F_t(z)$ with respect to 
the variables $z_1, z_2$.
It is easy to see that 
\BQ
K(0)&=& 1 \label{DE1} \\ \label{DE2}
\frac {\p }{\p t}K(t)&=&(A(z)+\triangledown A(z))K(t) 
\EQ

To show that $K(t)=H(t)$, it is enough to 
show that $H(t)$ also satisfies the equations (\ref{DE1}) and (\ref{DE2}) above. 
First when $t=0$, $F_t(z)=(z_1, z_2)$ and $H(0)={\mathcal J}(F)(z; 0)=1$.
So it only remains to check ({\ref{DE2}) for $H(t)$.

\BQ
\frac {\p }{\p t}H(t)&=&
\frac {\p }{\p t} \left | 
\begin{matrix} {\frac {\p F_1(z; t)}{\p z_1}}, &{\frac {\p F_1(z; t)}{\p z_2}} \\
{\frac {\p F_2(z; t)}{\p z_1}}, &{\frac {\p F_2(z; t)}{\p z_2}} \end{matrix} \right |
 \nno \\ \label{PE3}  
&=& \left |
 \begin{matrix} 
{\frac {\p^2 F_1(z; t)}{\p z_1\p t}}, 
&{\frac {\p F_1(z; t)}{\p z_2}} \\
{\frac {\p^2 F_2(z; t)}{\p z_1 \p t}}, 
&{\frac {\p F_2(z; t)}{\p z_2}} \end{matrix} \right |
+\left |
\begin{matrix} {\frac {\p F_1(z; t)}{\p z_1}}, &{\frac {\p^2 F_1(z; t)}{\p z_2\p t}} \\
{\frac {\p F_2(z; t)}{\p z_1}}, 
&{\frac {\p^2 F_2(z; t)}{\p z_2\p t}} \end{matrix} \right | 
\EQ

By Lemma \ref{L2.6}, we calculate the first term of (\ref{PE3}) as follows.

\BQ
&{}& \left |
 \begin{matrix} {\frac {\p^2 F_1(z; t)}{\p z_1\p t}}, 
&{\frac {\p F_1(z; t)}{\p z_2}} \\ \nno
{\frac {\p^2 F_2(z; t)}{\p z_1\p t}},
 &{\frac {\p F_2(z; t)}{\p z_2}} \end{matrix} \right |\\ \nno
&=& \left |
 \begin{matrix} {\frac {\p}{\p z_1}(a_1(z)\frac {\p}{\p z_1}+a_2(z)\frac {\p}{\p z_2})
F_1(z; t)}, &{\frac {\p F_1(z; t)}{\p z_2}} \\ \nno
{\frac {\p}{\p z_1}(a_1(z)\frac {\p}{\p z_1}+a_2(z)\frac {\p}{\p z_2})
F_2(z; t)}, 
&{\frac {\p F_2(z; t)}{\p z_2}} \end{matrix} \right | \nno\\ \nno
&=& \left |
 \begin{matrix}{A\frac {\p F_1(z; t)}{\p z_1}}, &{\frac {\p F_1(z; t)}{\p z_2}}\\
{A\frac {\p F_2(z; t)}{\p z_1}}, &{\frac {\p F_2(z; t)}{\p z_2}} \end{matrix} \right |
+ \left |
\begin{matrix}{\frac {\p a_1}{\p z_1}\frac {\p F_2(z; t)}{\p z_1}+
\frac {\p a_2}{\p z_1}\frac {\p F_2(z; t)}{\p z_2}}, 
&{\frac {\p F_2(z; t)}{\p z_2}}\\
{\frac {\p a_1}{\p z_1}\frac {\p F_2(z; t)}{\p z_1}+
\frac {\p a_2}{\p z_1}\frac {\p F_2(z; t)}{\p z_2}},
 &{\frac {\p F_2(z; t)}{\p z_2}} \end{matrix} \right | \nno \\  \label{PE4}
&=& \left |
 \begin{matrix}{A\frac {\p F_1(z; t)}{\p z_1}}, 
&{\frac {\p F_1(z; t)}{\p z_2}} \\ 
{A\frac {\p F_2(z; t)}{\p z_1}}, &{\frac {\p F_2(z; t)}{\p z_2}} \end{matrix} \right |
+\left ( \frac {\p a_1}{\p z_1}\right ) {\mathcal J} (F_t)
\EQ

Similarly, for the second term of (\ref{PE3}), we have
\BQ \label{PPE5} 
\left |
\begin{matrix} 
{\frac {\p F_1(z; t)}{\p z_1}}, &{\frac {\p^2 F_1(z; t)}{\p z_2\p t}} \\
{\frac {\p F_2(z; t)}{\p z_1}}, 
&{\frac {\p^2 F_2(z; t)}{\p z_2\p t}} 
\end{matrix} \right | 
= \left |
 \begin{matrix}
{\frac {\p F_1(z; t)}{\p z_1}}, &{A\frac {\p F_1(z; t)}{\p z_2}} \\
{\frac {\p F_2(z; t)}{\p z_1}}, &{A\frac {\p F_2(z; t)}{\p z_2}} 
\end{matrix} \right |
+ \left ( \frac {\p a_2}{\p z_2}\right ) {\mathcal J} (F_t)
\EQ

Combining (\ref{PE4}) and (\ref{PPE5}), we get
\BQ
\frac {\p }{\p t}H(t)&=&
A \left | 
\begin{matrix} 
{\frac {\p F_1(z; t)}{\p z_1}}, &{\frac {\p F_1(z; t)}{\p z_2}} \\
{\frac {\p F_2(z; t)}{\p z_1}}, &{\frac {\p F_2(z; t)}{\p z_2}} 
\end{matrix} \right |
+ (\frac {\p a_1}{\p z_1}+\frac {\p a_2}{\p z_2})\left | 
\begin{matrix} 
{\frac {\p F_1(z; t)}{\p z_1}}, &{\frac {\p F_1(z; t)}{\p z_2}} \\
{\frac {\p F_2(z; t)}{\p z_1}}, &{\frac {\p F_2(z; t)}{\p z_2} } 
\end{matrix} \right | \nno \\ \label{E2.32}
&=& (A+\triangledown A){\mathcal J}(F_t) 
\EQ

$(b)$ By formula \ref{MF} and Lemma \ref{L2.6}, it is easy to check that 
both sides of (\ref{G-MF}) satisfy equations (\ref{DE1}) and (\ref{DE2}).
\epfv

By the similar idea, we also can  get an exponential formulas for the Jacobian 
matrix $JF(t, z)$ of $F(t, z)$. First we fix the following notations: 
Let $Ja(z)$ be the Jacobian matrix of the n-tuple 
$(a_1(z), \cdots , a_n(z))$. Let $R_{Ja}$ be the operator over 
the algebra $M_{n\times n}(\bC [[z]])$, i.e. the $n \times n$ matrices 
with entries lying in $\bC [[z]]$, defined by multiplifying 
the matrix $Ja(z)$ from the right-hand side. In the following theorem, we also
view the differential operator $A(z)=a(z)\frac{\p}{\p z}$ 
as a differential operator of the algebra $M_{n\times n}(\bC [[z]])$, which acts
on the  matrices entry-wisely.

\begin{theo}\label{JM-MT}
 For any $U(z) \in M_{n\times n}(\bC [[z]])$, we have
\BQ\label{GJM-formula}
exp(tA+tR_{Ja}) U=U(F_t(z))JF_t(z)
\EQ
In particular, when $U$ is chosen to the identity matrix $Id$, we get
\BQ \label{JM-formula}
JF_t(z)=exp ( tA+ tR_{Ja} )\cdot Id
\EQ
\end{theo}

\pf  For any $1\leq i, j \leq n$, consider 
\BQn
\frac {\p}{\p t} \frac {\p F_i(t, z)}{\p z_j} 
&=&
\frac{\p }{\p z_j} \frac {\p F_i(t, z)}{\p t}  \\
&=& \frac{\p}{\p z_j}\sum_{k=1}^n a_k(z)\frac {\p F_i(t, z)}{\p z_k}  \\
&=&\sum_{k=1}^n \frac {\p a_k}{\p z_j} \frac {\p F_i(t, z)}{\p z_k} +
 (\sum_{k=1}^n a_k \frac {\p}{\p z_k} )\frac {\p F_i(t, z)}{\p z_j} \\
&=& \sum_{k=1}^n \frac {\p F_i(t, z)}{\p z_k} \frac {\p a_k}{\p z_j}
+ A \frac {\p F_i}{\p z_j} \\
\EQn

Hence we have
\BQ
\frac {\p}{\p t} JF_t(z)= (A+R_{Ja}) JF_t(z)
\EQ
By Lemma \ref{L2.6}, we also have $\frac {\p}{\p t} U(F_t(z))= A U(F_t(z))$. So
it is easy to see that the right hand side of (\ref{GJM-formula}) 
satisfies the equations
\BQ
\frac {\p}{\p t}\left ( U(F_t(z)) JF_t(z) \right ) &=&
  (A+R_{Ja})( U(F_t(z))JF_t(z)) \\
U(F_0(z))JF_0(z)=Id
\EQ
Hence (\ref{GJM-formula}) holds.
\epfv

\begin{rema}\label{rmk=1}
$(a)$ Note that the proofs of Theorem \ref{MT} 
and Theorem \ref{JM-MT} only need the condition $o(a(z))\geq 1$ instead of 
$o(a(z))\geq 2$. So
for any $A(z)=a(z)\frac {\p}{\p z}$ with $o(a(z))\geq 1$, set 
$F(t; z)=e^{tA(z)}z$, then the formulas in Theorem \ref{MT} 
and Theorem \ref{JM-MT} still hold. 

$(b)$ In particular, 
over the complex field $\bC$, it is straightforward 
to check that
$F(z)=e^{A(z)}z$ is a well defined formal power series 
and we can replace $t$ by $1$ in all the formulas in Theorem \ref{MT} 
and Theorem \ref{JM-MT}.
\end{rema}

Next we will derive a little bit more information about ${\mathcal J}F_t$.

\begin{propo} \label{lePE19}
\BQ \label{PE19}
\frac {\p}{\p t}{\mathcal J}(F_t)=
(A+(\triangledown a)(z)){\mathcal J}(F_t)=
(\triangledown a)(F_t){\mathcal J}(F_t)
\EQ
In particular, 
\BQ
A {\mathcal J}(F_t)=((\triangledown a)(F_t)-(\triangledown a)(z)){\mathcal J}(F_t)
\EQ
\end{propo}

\pf From (\ref{E2.32}), we see that $\frac {\p}{\p t}{\mathcal J}(F_t)=
(A+(\triangledown a)(z)){\mathcal J}(F_t)$. Let 
$L(z; t)=(A+(\triangledown a)(z)){\mathcal J}(F_t)$ and 
$R(z; t)=(\triangledown a)(F_t){\mathcal J}(F_t)$. Then by (\ref{MF}) and 
Lemma \ref{L2.6}, it is easy to see that
\BQ
\frac {\p L(z; t)}{\p t}&=&(A+(\triangledown a)(z))L(z; t) \\
\frac {\p R(z; t)}{\p t}&=&(A+(\triangledown a)(z))R(z; t)
\EQ
While $L(z; 0)=(\triangledown a)(z)=R(z; 0)$, Hence we must have
$L(z; t)=R(z; t)$.
\epfv

As an application of Theorem \ref{MT}, we give a new proof to the 
following result, 
which was first proved by M. Pittaluga in \cite{P} by using the theory of
formal Lie groups and Lie algebras.

\begin{corol}\label{C2.10}
${\mathcal J}(F)\equiv 1$ if and only if $\triangledown A\equiv 0$.
\end{corol}
\pf  First from (\ref{MF}), 
it is easy to see that if $\triangledown A\equiv 0$, then 
${\mathcal J}(F)\equiv 1$. Conversely, suppose that ${\mathcal J}(F)\equiv 1$.
Observe that $a(z)\in \bC_2[[z]]$, or in other words, the least degree of $a_i$ 
are at least $2$, therefore the operators $A=a(z)\frac {\p}{\p z}$ and
$A+\triangledown A$ increase the degree at least by one. 
If $\triangledown a(z)\neq 0$, say its lowest degree is $m$. Let $M$ be 
it the homogeneous part of degree $m$. From (\ref{MF}) for $t=1$, we have

\BQ
1&\equiv &{\mathcal J}(F)=e^{(a(z)
\frac {\p}{\p z}+\triangledown a(z))} \cdot 1 \nno \\
&=&1+(a(z)\frac {\p}{\p z}+\triangledown a(z))\cdot 1+
\sum_{i\geq 2}\frac 1{k!}(a(z)\frac {\p}{\p z}+\triangledown a(z))^{k-1}
\triangledown a(z) \nno \\ \label{PE6}
&=& 1+M+\text{high degree terms} 
\EQ

Clearly $M=0$, contradiction.

Another way to prove the result above is the following: 
Consider the ``deformation'' $F_t(z)$ of $F$ as before. Notice the Jacobian
${\mathcal J}(F_t)\in \bC [t][[z]]$ and ${\mathcal J}(F_t)={\mathcal J}(F^{[k]})$ when $t=k$,
 for any
$k\in \bN$. Now since
${\mathcal J}(F)(z, 1)\equiv 1$, 
then, by the chain rule,  ${\mathcal J}(F^{[k]})(z)\equiv 1$. This implies that 
${\mathcal J}(F_t)\equiv 1$,  when $t=k$ for any $k\in \bN$. Hence
${\mathcal J}(F_t)$ itself must be identically $1$, for as a polynomial of $t$, the
coefficient of any monomial of positive degree of $F$ can not have infinitely roots
unless it is zero.
In particular, ${\mathcal J}(F_t)$ does not depends on t. So we have

\BQ
0=\frac {\p}{\p t} \lbar_{t=0} {\mathcal J}F(z; t)=
(a^{(z)}\frac {\p}{\p z}+\triangledown a(z)){\mathcal J}F(z; 0)=\triangledown a(z)
\EQ

\epfv

From the arguments in the proof for the Corollary above,  
or by the Corollary itself, we have
 
\begin{corol}
For any $F \in {\mathcal F}_1$, if ${\mathcal J}(F)\equiv 1$. Then
${\mathcal J}(F_t)\equiv 1$.
\end{corol}

\renewcommand{\theequation}{\thesection.\arabic{equation}}
\renewcommand{\therema}{\thesection.\arabic{rema}}
\setcounter{equation}{0}
\setcounter{rema}{0}

\section{Some Explanations and Applications} \label{Sec-3}

At the first glance, the formulas we proved in Theorem \ref{MT}
and Theorem \ref{JM-MT} are a little mysterious. 
Here we try to give a little 
explainations to these two formulas.

First the exponenitial formula (\ref{MF}) reminds us the following 
so called Liouville's formula in linear algebra. Namely, for any 
$n\times n$ matrix $M\in M_{n\times n}(\bC)$, then
\BQ \label{LAF}
Det \, e^M=e^{Tr M}
\EQ

Actually we will see that the formula (\ref{MF}) can be viewed as 
a generalization of the Liouville's formula above.

First we define the embedding 
\BQ
\Phi : M_{n\times n}(\bC) &\to & {\mathcal D}(z) \\
M=(m_{ij}) &\to & \sum_{i, j=1}^n m_{ij}z_i \frac {\p}{\p z_j}
\EQ
where ${\mathcal D}(z)$ is the Lie algebra of 
the derivations of $\bC [[z]]$.
It is very easy to check that the linear map
$\Phi : M_{n \times n} \to {\mathcal D}(z)$ is an injective 
homomorphism of Lie algebras.

\begin{lemma} Let $F(z)= exp(\Phi (M))z$. Then

$(a)$ $J(F)= e^M$.

$(b)$ $F(z)=e^M z$.

$(c)$  ${\mathcal J}(F)=e^{Tr M}$.
\end{lemma}

\pf 
Note that $J \Phi (M)=M$ and 
$\triangledown \Phi (M)=Tr M$. By Remark \ref{rmk=1},
we can apply formula (\ref{JM-formula})
to the map $F$, we get
\BQn
J(F)&=& e^{\Phi(M)+ R_{J\Phi(M)}}I_{n \times n} \\
&=& e^{R_{M}}e^{\Phi(M)} I_{n \times n} \\
&=& e^{R_{M}}I_{n \times n} \\
&=& e^M
\EQn
where the second equality above follows from the fact that the operators 
$\Phi(M)$ and $R_{J\Phi(M)}$ commutes with each other. So we have proved $(a)$. 
$(b)$  follows immediately from $(a)$.

To prove $(c)$, we apply formula (\ref{MF}) to $F$, we get
\BQn 
{\mathcal J}(F)&=& e^{\Phi(M)+ \triangledown \Phi (M)} \cdot 1 \\
&=& e^{\triangledown \Phi (M)}e^{\Phi(M)}\cdot 1 \\
&=& e^{Tr (M)}
\EQn
\epfv

Combine $(a)$ and $(c)$ in the lemma above, we recover formula (\ref{LAF}).
Therefore formula (\ref{MF}) and (\ref{JM-formula}) can be viewed 
as some generalizations of the Liouville's formula (\ref{LAF}).

One of the motivations for the present work is that we believe the exponential 
formulas (\ref{EF1}), (\ref{MF}) and Corollary \ref{C2.10} are closely related with 
the well known Jacobian Conjecture. In the rest of section, we will consider 
some applications to the Jacobian Conjecture.

From Proposition \ref{6P2.1}, Lemma \ref{L2.5} and Corollary \ref{C2.10},
we see that
the Jacobian Conjecture is equivalent to the following pure algebraic
problem.

\begin{conj}
Let $a(z) \in \bC_2 [[z_1, z_2, \cdots , z_n]]$ and $\triangledown a(z)=0$. Then
$F(z)=exp ({a(z)\Pz}) z\in (\bC [z])^n$ if and only if 
$G(z)=exp ({-a(z)\Pz})z       \in (\bC [z])^n$.
\end{conj}

In the case when $a(z)$ is even, we have a very simple answer to the conjecture above.

\begin{propo}\label{P-even}
$a)$ For any $F\in {\mathcal F}_1$, let $G$ be its formal inverse. Then
$G(z)=-F(-z)$ if and only if $a(z)$ is even.

$b)$ If $F$ satisfies the conditions in the Jacobian Conjecture and
$a(z)$ is even, then $G$ is also a polynomial map.
\end{propo}

\pf  Clearly $b)$ is an immediate consequence of $a)$. For $a)$, 
Suppose $a(z)$ is even, then,
replacing  $z$ by $-z$ in (\ref{EF1}), we get
\BQ
F(-z)&=& exp(a(-z)\frac {\p}{\p (-z)})(-z) \nno \\
&=& -exp(-a(z)\frac {\p}{\p z})z \nno \\ \label{P3.3.1}
&=& -G(z)
\EQ

Conversely, suppose the formal inverse $G(z)=-F(-z)$. 
Let $B=b(z)\Pz$ be the 
associated formal differential operator of $G$, i.e. 
\BQ\label{P3.3.2}
G(z)=exp(B(z))z 
\EQ
By the 
uniqueness of $B$, we have $B(z)=-A(z)$. On the other hand,
 from (\ref{P3.3.1}), we get
\BQ\label{P3.3.3}
-F(-z)&=& exp(A(-z))z
\EQ
Comparing (\ref{P3.3.2}) and (\ref{P3.3.3}), we have
$A(-z)=B(z)=-A(z)$. Therefore $a(z)$ must be even.
\epfv

As an immediate consequence, we have the following:
\begin{corol}
With the same notation above, if
 $a(z)$ is even and $F= e^A z$ are polynomials, then
$\triangledown a(z)=0$. 
\end{corol}

Note that this is not true for arbitrary formal power series $a(z)$.

Finally, we give a new proof for a theorem of Bass, Connell and Wright in \cite{BCW}.

\begin{theo}\label{Thm-BCW}  \cite{BCW}
Let $F(z)=z+H(z)$ be a polynomial map with $H(z)$ being homogeneous of degree $d\geq 2$. 
If $J(H)^2=0$, then the formal inverse map $G=z-H(z)$. 
\end{theo}

Note that $J(H)^2=0$ implies that ${\mathcal J}(F)=1$. Thus the Jacobian
Conjecture is true in this case.

\pf First note that $J(H)z=dH(z)$, since $H(z)$ is homogeneous of degree $d$. 
From $J(H)^2=0$, we have $0=J(H)^2z=d J(H)H$, hence $J(H)H=0$.

Now let $A(z)=a(z)\frac {\p}{\p z}$ be the  calculate the associated formal differential 
operator.  Write $a(z)=\sum_{k=2}^\infty a^{(k)}(z)$, where $a^{(k)}(z)$ is 
homogeneous of degree $k$. By incursive formula (\ref{Incursive}), it is 
easy to see that $a^{(k)}(z)=0$ if $k\neq m(d-1)+1$ for some $m>0$. For 
$k= m(d-1)+1$ with  $m>0$, we have $a^{(d)}(z)=H(z)$ and
\BQn
a^{(d)}(z) &=& H(z) \\
a^{(2d-1)}(z) &=& -\frac 12 (H(z)\frac {\p}{\p z})^2 z \\
&=& -\frac 12 (H(z)\frac {\p}{\p z})H(z)\\
&=& -\frac 12 JH(z) \cdot H(z)\\
&=& 0
\EQn
By Mathematical Induction and incursive formula (\ref{Incursive}), 
it is easy to show that 
$a^{(m(d-1)+1)}=-\frac 1{m!} (H(z)\frac {\p}{\p z})^mz=0$ for any $m\geq 2$. Therefore, 
we have $a(z)=H(z)$ and $A(z)=H(z)\frac {\p}{\p z}$. Note that
$A^2(z)=0$, so the formal inverse $G(z)$ of $F(z)$ is given by
$G(z)=e^{-A}z=z-H(z)$.
\epfv

\renewcommand{\theequation}{\thesection.\arabic{equation}}
\renewcommand{\therema}{\thesection.\arabic{rema}}
\setcounter{equation}{0}
\setcounter{rema}{0}

\section{Some Open Problems}\label{Sec-4}

For the case of two variables, by using the residue and intersection theory
in complex algebraic geometry,
 the author in \cite{Z2} shows that, to 
prove the Jacobian Conjecture,   
it will be enough to consider the following special polynomial maps 
$F\in {\mathcal F}_1$.

Let $r(x)$ be a monic polynomial of degree $N+1>1$ with distinct roots and
$\lambda (x)$ and $\mu (x)$ are unique polynomials satisfying 

$a)$
\BQ
r(x)\mu (x)+r'(x)\lambda (x) &=& 1
\EQ

$b)$ $\deg \mu (x)\leq N-1$ and $\deg \lambda (x)\leq N$. 

Consider $F=(F_1, F_2)$,  where
\BQ
F_1(z_1, z_2)&=&r(z_1)H_1(z_1, z_2)-z_2\lambda (z_1)K_2(z_1, z_2) \\
F_2(z_1, z_2)&=&r(z_1)H_2(z_1, z_2)+z_2\lambda (z_1)K_1(z_1, z_2) 
\EQ
where $H_i$ and $K_i$ are polynomials in $z=(z_1, z_2)$ and 
satisfy $H_1K_1+H_2K_2=1$. Furthermore, 
without losing any generality, we also can assume that $F\in {\mathcal F}_1$.
Then the Jacobian Conjecture is equivalent to the following

\begin{conj}
Let $F=(F_1, F_2)$ as above, $A=a(z)\Pz$ be the associated 
formal differential operator
of $F$, then $\triangledown A\neq 0$.
\end{conj} 

Finally we ask the following very important and
 interesting question.
(This question for the case of one variable was first suggested to 
the author by Y-Z. Huang), namely, if
the analytic map $F$ is well defined in an open 
neighborhood of $0\in \bC^n$, is $a(z)$ convergent near $0\in \bC^n$ ?

This is unknown both in the case of 
one variable and in the case $F$ is a
polynomial map with ${\mathcal J}(F)\equiv 1$. We believe the following conjecture
is true, but we do not have much evidence.

\begin{conj}
If $F$ is convergent near $0\in \bC^n$, then so is $a(z)$.
\end{conj}

The converse of the conjecture above is very easy to prove.

\begin{propo}
Suppose $a(z)\in \bC_2 [[z]]$ is convergent near point $0\in \bC^n$, 
then so is the formal power series $F(z)=e^{a(z)\Pz}z$.
\end{propo}
\pf  
Consider the deformation $F(z; t)=e^{ta(z)\Pz}z$, which satisfies the following
differential equations

\BQ \label{Deq}
\frac {\p}{\p t}F(z; t) &=& a(z)\Pz F(z; t)\\
F(z; 0) &=& z \label{Deq-Con}
\EQ

It is well known in PDE that the differential equation (\ref{Deq}) with
condition (\ref{Deq-Con}) has a unique
analytic solution. Then as a formal power series solution of (\ref{Deq}), $F$ is
convergent near $0\in \bC^n$.
\epfv

{\small \sc Department of Mathematics, Washington University in St. Louis,
St. Louis, MO 63130-4899}

{\em E-mail}: zhao@math.wustl.edu

\end{document}